\documentclass[12pt]{amsart}

\usepackage{amssymb}

\usepackage{enumitem}

\usepackage{graphicx}
\usepackage{tikz}
\usetikzlibrary{positioning}
\usepackage{array}
\usepackage{geometry}

\makeatletter
\@namedef{subjclassname@2020}{%
  \textup{2020} Mathematics Subject Classification}
\makeatother

\usepackage[T1]{fontenc}

\newtheorem{theorem}{Theorem}
\newtheorem{corollary}[theorem]{Corollary}
\newtheorem{lemma}[theorem]{Lemma}
\newtheorem{proposition}{Proposition}
\newtheorem{example}{Example}
\newtheorem{remark}{Remark}
\newtheorem{fact}{Fact}
\newtheorem{definition}{Definition}

\theoremstyle{definition}

\numberwithin{equation}{section}

\begin{document}

\title{Quotient structure of ultrafilters}
\author{Manoranjan Singha and Rohan Pradhan}


\begin{abstract}
 Various notions of divisibility of ultrafilter by ultrafilter are flourished by elegant minds during last decade leading to a good literature on divisibility relations on ultrafilters. This article extends the study by shaping  the quotient structure resulting particularly from strong divisibility on $\beta\mathbb N$ that provides a natural analogue of quotient formation in $\mathbb N$. Basic properties of quotient structure are established, and it is showed that the divisibility profile of an ultrafilter $u$ is preserved through the products $\frac{u}{v}\cdot v$ and $v\cdot \frac{u}{v}$, where $v$ is a divisor  of $u$. A characterization of self-divisible ultrafilters is obtained in terms of quotient structure arising from the division of underlying ultrafilter by any of its divisors. For self-divisible multiplicative idempotent $u$, it is proved that self-divisibility is inherited by the quotient $\frac{u}{u}$. For self-divisible non-principal ultrafilter $u$, a necessary and sufficient condition for strong divisibility of the quotient $\frac{u}{u}$ by any prime ultrafilter is established using exponentiations of prime ultrafilter. Multiplicative $\Delta$-sets and $\Delta^{*}$-sets are characterized via quotients of self-divisible ultrafilters. Finally, multiplicative homomorphisms are characterized through the compatibility of the quotient structure with Stone–Čech extensions.
\end{abstract}




\maketitle

\section{Introduction}
The Stone-
ech Compactification of the discrete space $\mathbb{N} $ denoted by $ \beta \mathbb{N}$ is the space of all the ultrafilters on $ \mathbb{N} $ endowed with the Stone topology generated by the basic open sets of the form $ \widehat{A} = \{ u \in \beta \mathbb{N}: A \in u\} $. Each $ n \in \mathbb{N}$ is identified with the ultrafilter $ u_n = \{ A \subset \mathbb{N} : n \in A \}  $ and members of $ \beta \mathbb{N}\setminus \mathbb{N}$ are called non-principal ultrafilters. Every member of non-principal ultrafilter is infinite.
The usual addition and  multiplication on $\mathbb{N}$ can be extended to $\beta \mathbb{N}$;
\\ For $u,v \in \beta \mathbb{N},  u +v = \{A \subset N : \{x \in \mathbb{N}: A-x \in v \} \in u\}$, where $A -x = \{y \in \mathbb{N}: y+x \in A \}$ , $ u\cdot v =  \{A \subset N : \{x \in \mathbb{N}: A/x \in v \} \in u\}$ and $A /x = \{y \in \mathbb{N}: y\cdot x \in A \}$ . It is well-established that with these operations,  $ \beta \mathbb{N}$ forms a compact Hausdorff right topological semigroup \cite{4}. Since every compact Hausdorff right topological semigroup contains an idempotent \cite{4}, the semigroup $ \beta \mathbb{N}$ admits both additive and multiplicative idempotents which play important role in this article. Though $\beta \mathbb{N}$ contains both additive and multiplicative idempotents, there does not exist any ultrafilter in $\beta \mathbb{N}$ which is additive as well as multiplicative idempotent \cite{4}. By uniquely extending standard arithmetic operations like addition and multiplication from $\mathbb{N}$ to $\beta \mathbb{N}$, the algebraic and combinatorial properties of ultrafilters can be analyzed. Several recent works have significantly contributed to the development of ultrafilter theory and its interaction with topological dynamics and combinatorial structures \cite{5}, \cite{6}. In recent years, the study of divisibility of ultrafilters has emerged as one of the most active and intriguing areas of research in the theory of $\beta \mathbb{N}$. Some definitions and results / facts are given below as ready references which are used in the sequel.
\begin{fact} \label{fact 1}\cite{4}
	Every map $ f : \mathbb{N} \rightarrow \mathbb{N} $ extends uniquely to a continuous map $ \widetilde{f} : \beta \mathbb{N} \rightarrow \beta \mathbb{N} $ and for every $ u  \in \beta \mathbb{N} ,    \widetilde{f}(u) = \{ A \subset \mathbb{N} : f^{-1}(A) \in u\}$.
\end{fact}
In view of Fact \ref{fact 1}, B {S}obot introduced the relation $\widetilde{\mid}$ as an extension of standard divisibility $\mid $  on $ \mathbb{N}$ to $\beta \mathbb{N}$ in \cite{7}. More about this relation $\widetilde{\mid}$  has been studied in \cite{7}, \cite{8}, \cite{9}, \cite{10} and \cite{11}. In \cite{10}, B \v{S}obot introduced a new type of divisibility relation $\mid^{s}$ as a stronger version of $\widetilde{\mid}$-divisibility namely strong divisibility and showed that this divisibility relation is compatible with the congruence relation.
\begin{definition}\label{definition 1}
	\cite{3}	For $u \in \beta \mathbb{N}, \ D(u) = \{x \in \mathbb{N} : x\mathbb{N} \in u\}$.
\end{definition}
\begin{definition}\label{definition 2}
	\cite{7}	For $u,v \in \beta \mathbb{N},\ v \widetilde{\mid} u$ if for every $ A \in v, \ \bigcup_{n \in A} n\mathbb{N} \in u$.
\end{definition}
\begin{definition}\label{definition 3}
	\cite{7}	For $u,v \in \beta \mathbb{N}, u$ is said to be right divisible by $v$ if there exists $ w \in \beta \mathbb{N}$ such that $ u = vw$. Similarly,  $u$ is said to be left divisible by $v$ if there exists $ w \in \beta \mathbb{N}$ such that $ u = wv$. Both right and left divisibilities imply $\widetilde{\mid}$ - divisibility but the converse does not hold in general.
\end{definition}
\begin{definition}\label{definition 4}
	\cite{3} 	Given $u, v \in \beta \mathbb{N}, u $ is strongly divisible by $v$  denoted as $ v \mid^{s} u $ whenever $ D(u) \in v$.\\ It is observed that for any $ n \in \mathbb{N}$ and $u \in \beta \mathbb{N}, n \mid^{s} u $  iff there exists $ A \subset \mathbb{N}$ such that $ A \in u $  and $ n \mid m, \forall m \in A $. Replacing $n$ by an arbitrary ultrafilter $v$, it follows that $ v \mid^{s} u $ iff there exists $ B \subset \mathbb{N}$ such that $ B \in v $ and for every $ b \in B, \ \exists A \in u $ such that $ b \mid a, \ \forall a \in A$.
	\\Note: For principal ultrafilters $ \mid^{s}$ and $ \widetilde{\mid} $ are equivalent. Also for any $u \in \beta \mathbb{N} $ and $n \in \mathbb{N}, n \mid^{s} u \iff n \widetilde{\mid} u \iff n\mathbb{N} \in u $.   
\end{definition}
\begin{fact}\label{fact 2}
	\cite[Theorem 6.6]{10}	Let $u,v, w \in \beta \mathbb{N}$. Then (i) $ v \mid^{s} u \implies   v \mid^{s} u\cdot w $, (ii) $ v \mid^{s} u \implies w\cdot u$.
\end{fact}
\begin{definition}\label{definition 5}
	\cite{3} An element $ u \in \beta \mathbb{N}$ is said to be self-divisible if $ u \mid^{s} u$, that is, if $ D(u) \in u$.
\end{definition}

\begin{definition}\label{definition 6}
	\cite{9} An ultrafilter $u$ is called a prime ultrafilter if $ P \in u $ where $P = \{x \in \mathbb{N}: x \text{ is prime} \} $ .
	
\end{definition}
Throughout the paper P denotes the set of prime numbers.

\begin{definition}\label{definition 7}
	\cite{8} The class of ultrafilters divisible by all the members of $\beta\mathbb{N}$ is denoted by $MAX$ and is defined as $MAX = \bigcap_{n \in \mathbb{N}} \widehat{n\mathbb{N}}$. In terms of $ D(u), MAX  = \{u \in \beta \mathbb{N}: D(u) = \mathbb{N}\}$.
	\\ Note : $ MAX$ is a non-empty closed subset of $\beta \mathbb{N}$, in fact every additive idempotent is a member of $ MAX$ \cite{3}. 
\end{definition}
\begin{definition}\label{definition 8}
	\cite{4}Let $A \subset\mathbb{N}$. $A$ is said to be an additive $IP$-set if there exists a sequence $(x_n)$ in $ \mathbb{N} $ such that $FS(x_n) \subset A $, where $ FS(x_n) = \{\sum_{x \in F} x: F$ is a finite non-empty subset of $\mathbb{N}\} $. Definition of multiplicative $IP$-set is similar, one needs to replace finite sum by finite product.
\end{definition}
\begin{fact} \label{fact 3}
	\cite{4}Any subset $A$ of $\mathbb{N}$ is an additive $IP$-set iff $A \in u$, for some additive idempotent ultrafilter $u$. [Analogous result holds for multiplicative $IP$-set]
\end{fact}
\begin{definition}\label{definition 9}\cite{3}
	Given $ u \in \beta \mathbb{N}$, $ \varphi_u : P \longrightarrow \omega +1$  is defined as the function sending p to max$\{k \in \omega : p^k\mathbb{N} \in u\}$, if this exists and to $ \omega$ otherwise.
\end{definition}
Following the convention \ $p^{\omega+1} \mathbb{Z} = \{0\}$ \cite{3}, $ p^{\omega+1} \mathbb{N} = \phi $ is considered throughout the article. This leads to define $ D(u) $ in terms of the function $\varphi_u$ as $ D(u) = \bigcap_{p\in P } (p^{\varphi_u(p)+1}\mathbb{N})^c$ and $ (D(u))^c = \bigcup_{p \in P} p^{\varphi_u(p)+1}\mathbb{N}$.
\begin{definition}\label{definition 10}\cite{3}
	For $ u \in \beta \mathbb{N}$, the map $\varphi_u$ is  said to be finite if $ \varphi_u^{-1}(\{\omega\}) = \phi$ and $\varphi_u^{-1}(\mathbb{N}) $ is finite.  Likewise, the map $\varphi_u$  is said to be cofinite if $ \varphi_u^{-1}(\{0\}) \cup \varphi_u^{-1}(\mathbb{N}) $ is finite, that is,  $ \varphi_u^{-1}(\{\omega\})$ is cofinite in P. 
\end{definition}

\section{Focused area of the article }
So far as in the existing literature there is no concept of quotient though various type of divisibilities of ultrafilters by ultrafilters are there in \cite{3}, \cite{7}, \cite{8}, \cite{9}, \cite{10}, \cite{11}. In this article, the stronger version $ \mid^{s} $ of $ \widetilde{\mid}$-divisibility is adopted  which provides a framework to study self-divisibility due to the lack of reflexivity of $ \mid^{s} $ though the main attention is given to explore the quotient structure in either cases when an ultrafilter $u$ is divisible or not divisible by an ultrafilter $v$.
\\ For $ m,n \in \mathbb{N} , n \mid m \iff m \in n\mathbb{N} \iff \{x \in \mathbb{N}: x\mathbb{N}\in u_m\} \in u_n$, where $ u_m $ and $ u_n $ denote the principal ultrafilters corresponding to $ m$ and $n$ respectively. This observation motivates to study quotient-type structure on $\beta \mathbb{N} $.
\begin{definition}\label{definition 11}
	For $ u,v \in \beta \mathbb{N}$, define $ \frac{u}{v} = \{A \subset \mathbb{N}:\{x \in \mathbb{N}:  xA \in u\}\in v \}$.
	
\end{definition}

\begin{remark}\label{remark 1}
	
	For $ u, v \in \beta \mathbb{N},\frac{u}{v} \in \beta \mathbb{N} \iff \frac{u}{v} \neq \phi. $

\end{remark}
From Definition \ref{definition 4}, whenever $ v\mid^{s}u, D(u) \in v$. By the definition of  $\frac{u}{v}, D(u) =\{x \in \mathbb{N}: x\mathbb{N} \in u \} \in v $ implies $ \mathbb{N} \in \frac{u}{v}$. This allows us to define the stronger divisibility $\mid^{s}$ in terms of $\frac{u}{v}$ instead of $D(u)$ \cite{3}.

\begin{definition}\label{definition 12}
	Given $ u,v \in \beta \mathbb{N}, u $ is said to be strongly divisible by $v$ if   $\frac{u}{v}$ is non-empty. So whenever $u$ is not strongly divisible by $v$, the quotient $\frac{u}{v}$ becomes empty.
	
\end{definition}
\begin{fact}\label{fact 4}
	For $ u \in \beta\mathbb{N} \setminus \mathbb{N}$ and $ v \in \beta\mathbb{N}$ such that $ v \mid^{s} u , \ \frac{u}{v} \in \beta \mathbb{N}\setminus \mathbb{N}$.
\end{fact}
\begin{fact}\label{fact 5}
	For any $u,v \in \beta \mathbb{N}$ with $v \mid^{s} u$ , $D(v) \subset D(u)$ and $D(\frac{u}{v}) \subset D(u)$. 
\end{fact}
\begin{proposition}\label{proposition 1}
	Let $u \in \beta \mathbb{N}$. Then for any $v \in \beta \mathbb{N}$ with $v \mid^{s} u$, $ u \in MAX$ if and only if $\frac{u}{v} \in MAX $.
	\begin{proof}
		Take any $v \in \beta \mathbb{N}$ and suppose that $\frac{u}{v} \notin MAX $. Then there exists $k \in \mathbb{N}$ such that $k\mathbb{N} \notin \frac{u}{v}$. This further implies $\{x \in \mathbb{N}: x(\mathbb{N}\setminus k\mathbb{N})\in u\} \in v$. Choose $x \in \mathbb{N}$ such that $x(\mathbb{N}\setminus k\mathbb{N})\in u$. Since $x(\mathbb{N}\setminus k\mathbb{N}) \subset \mathbb{N}\setminus xk\mathbb{N}, \ \mathbb{N}\setminus xk\mathbb{N} \in u$. But this contradicts the fact that $u \in MAX$. Hence $\frac{u}{v} \in MAX$.
		\\Using Fact \ref{fact 5}, the converse follows immediately.
	\end{proof}
\end{proposition}
\begin{fact}\label{fact 6}
	Every additive idempotent ultrafilter belongs to $MAX$. Consequently, if $u$ is additive idempotent, then for any $v \in \beta \mathbb{N}$,  $\frac{u}{v}$ also belongs to $MAX$. Thus the quotient obtained on division of additive idempotent by any ultrafilter exhibits only trivial behavior.
\end{fact}
\begin{fact}\label{fact 7}\cite{3}
	For any $u \in \beta \mathbb{N}$,  the set of divisors of  $u$  in $\beta \mathbb{N}$ coincides with the set $\widehat{D(u)}$.
\end{fact}

For natural numbers $m,n $, whenever $ n \mid m $, there exists a quotient $ q = \frac{m}{n} \in \mathbb{N} $ such that $ m = \frac{m}{n}  \cdot n =  n \cdot \frac{m}{n} $. Consequently, the divisibility structure of $m$ is completely  reflected through that of products $  \frac{m}{n}  \cdot n$ and $  n \cdot \frac{m}{n} $. This naturally leads to a broader question : if both $m$ and $n$ are replaced by arbitrary ultrafilters, does the quotient structure still preserves the divisibility structure in a similar way ? 
\\ The following theorem shows that the answer is indeed affirmative.
\begin{theorem}\label{theorem 1}

	Let $u,v \in \beta \mathbb{N}$ be such that $ v \mid^{s} u$. Then for any $ w \in \beta \mathbb{N}$, the following equivalences hold \\$ w \mid^{s} u \iff w \mid^{s} \frac{u}{v} \cdot v \iff w \mid^{s} v \cdot \frac{u}{v}$. 
	\begin{proof}   From Fact \ref{fact 7} it follows that in order to establish the given equivalences it is sufficient to show that $ D(u) = D(\frac{u}{v} \cdot v) = D( v \cdot \frac{u}{v}) $. The proof is divided into four parts.
		\\ \textbf{1$\textsuperscript{st} $ part} : $  D(\frac{u}{v} \cdot v) \subset D(u)$. \\If possible let there exists  $x = p_1^{k_1} \cdot p_2^{k_2} \cdot \cdot \cdot p_n^{k_n} $, where $ p_1, p_2, ..., p_n $ are distinct primes and $k_i(1 \leq i \leq n) \in \mathbb{N}$, such that $ x \in  D(\frac{u}{v} \cdot v)$ but $ x \notin D(u) $. So there exists an $ i \in \{1,2,...,n\} $   for which $ p_i ^{k_i} \notin D(u) $. Setting  $ p_i = p $ and $ k_i = k $ we have  $ p^{k} \notin D(u)$. Consequently $\varphi_u(p) < k$ and $ 0 \leq \varphi_\frac{u}{v}(p) <k $. Let $\varphi_\frac{u}{v}(p) = m $. If $ m = 0 $ then $ p\mathbb{N} \notin \frac{u}{v}$ which means $\mathbb{N}\setminus p\mathbb{N} \in \frac{u}{v}$. Since $ x \in  D(\frac{u}{v} \cdot v), \ p^{k}\mathbb{N} \in \frac{u}{v} \cdot v $. Then $ (\mathbb{N}\setminus p\mathbb{N}) \cap \{y \in \mathbb{N} : (p^{k}\mathbb{N})/y \in v \} \in \frac{u}{v}$. Let $ y \in  (\mathbb{N}\setminus p\mathbb{N}) \cap \{y \in \mathbb{N} : (p^{k}\mathbb{N})/y \in v \}$. 
		\\ $ p \nmid y \implies (p^{k}\mathbb{N})/y = p^{k}\mathbb{N}$  but $ p^{k}\mathbb{N} \notin v $. Therefore $(\mathbb{N}\setminus p\mathbb{N}) \cap \{y \in \mathbb{N} : (p^{k}\mathbb{N})/y \in v \} = \phi $, a contradiction. 
		\\Now let us consider the case $ 0 < m < k$. Then $ p^{m}\mathbb{N} \in \frac{u}{v}$ but $ p^{m+1}\mathbb{N} \notin \frac{u}{v}$.
		\\ $  p^{m}\mathbb{N} \in \frac{u}{v}, p^{m+1}\mathbb{N} \notin \frac{u}{v}$,  $  p^{k}\mathbb{N} \in \frac{u}{v} \cdot v \implies ( p^{m}\mathbb{N}\setminus p^{m+1}\mathbb{N})  \cap \{y \in \mathbb{N} : (p^{k}\mathbb{N})/y \in v \} \in \frac{u}{v} $. Let $ y \in ( p^{m}\mathbb{N}\setminus p^{m+1}\mathbb{N})  \cap \{y \in \mathbb{N} : (p^{k}\mathbb{N})/y \in v \} $. Then $ (p^{k}\mathbb{N})/y  = p^{k-m}\mathbb{N} \in v $. Furthermore, $(p^{k-m}\mathbb{N}) \cap \{ t \in \mathbb{N}: t(p^{m}\mathbb{N}) \in u \} \in v$. But $ t_{0} \in (p^{k-m}\mathbb{N}) \cap \{ t \in \mathbb{N}: t(p^{m}\mathbb{N}) \in u \} \implies p^{k}\mathbb{N} \in u $, which is not possible and thus $ (p^{k-m}\mathbb{N}) \cap \{ t \in \mathbb{N}: t(p^{m}\mathbb{N}) \in u \} = \phi $, a contradiction. 
		\\ Therefore it leads to conclude that $ D(\frac{u}{v} \cdot v) \subset D(u) $.
		\\\textbf{ 2$\textsuperscript{nd} $ part} - $ D(v \cdot \frac{u}{v}) \subset D(u) $. \\Proof is similar to the 1$\textsuperscript{st} $ part, one just needs to  switch the roles of $\frac{u}{v} $ and $v$.
		\\ \textbf{3$\textsuperscript{rd} $ part} - $ D(u) \subset D(v \cdot \frac{u}{v}) $.
		\\ Suppose there exists  $ x = p_1^{k_1} \cdot p_2^{k_2} \cdot \cdot \cdot p_n^{k_n} $, where $ p_1, p_2, ..., p_n $ are distinct primes and $k_1, k_2, ..., k_n \in \mathbb{N}$, such that $ x \in D(u) $ but $ x \notin D(v \cdot \frac{u}{v}) $. Then there exists a prime $ p_i , 1 \leq  i \leq n $ such that  $ p_i ^{k_i} \notin D(v \cdot \frac{u}{v}) $. Setting  $ p_i = p $ and $ k_i = k $, we have $p^{k}  \notin D(v \cdot \frac{u}{v}) $ and $p^{k} \in D(u)$. Neither $p^{k}\mathbb{N} \in v$  nor $p^{k}\mathbb{N} \in \frac{u}{v}$.
		\\ $p^{k}\mathbb{N} \notin v,  p^{k}\mathbb{N} \notin \frac{u}{v}$, $ v \mid^{s} u \implies (\mathbb{N} \setminus p^{k}\mathbb{N} ) \cap \{\ y \in \mathbb{N}:  y(\mathbb{N} \setminus p^{k}\mathbb{N} ) \in u \} \cap D(u) \in v$. Consequently, $( p^{k}\mathbb{N} + i ) \in v $, where $ 1 < i < p^{k} - 1 $ and $ p \mid i $. Since $ p \mid i $, there exists $ a \in \mathbb{N} $ such that $ 1\leq a < k, \ p^{a} \mid i $ but $p^{a+1} \nmid i $.
		Let $ k - a = b $. Then $p^{b}\mathbb{N} \notin \frac{u}{v}$. So there exists $ j \in \mathbb{N}$ such that $ p^{b}\mathbb{N} + j \in \frac{u}{v} $ and $ 1\leq j \leq p^{b} - 1$. Moreover $( p^{k}\mathbb{N} + i ) \cap \{y \in \mathbb{N}: y(p^{b}\mathbb{N} + j) \in u \} \in v $. Pick $ y \in ( p^{k}\mathbb{N} + i ) \cap \{y \in \mathbb{N}: y(p^{b}\mathbb{N} + j) \in u \}$. 
		Then $ y = p^{k}t + i $ for some $ t \in \mathbb{N}$ and $( p^{k}t + i)(p^{b}\mathbb{N} + j) \in u $. But $ ( p^{k}t + i)(p^{b}\mathbb{N} + j) \cap p^{k}\mathbb{N} = \phi (\because  p^{k} \nmid ij) $, which leads to the contradiction. Therefore it follows that $ D(u) \subset D(v \cdot \frac{u}{v}) $.
		\\ \textbf{4$\textsuperscript{th} $ part} - $ D(u) \subset D( \frac{u}{v} \cdot v) $.
		\\ The proof is similar to that of the  3$\textsuperscript{rd} $ part.
		\\ Combining all four parts , we have $ D(u) = D(\frac{u}{v} \cdot v) = D( v \cdot \frac{u}{v}) $. Hence the result.
	\end{proof}
	
\end{theorem}

\begin{corollary}\label{corollary 2}
	Let $u,v \in \beta \mathbb{N}$ be such that $ v \mid^{s} u$. Then $ u \mid^{s} \frac{u}{v} \cdot v$ if and only if $u$ is self-divisible. (Similar result holds on replacing $\frac{u}{v} \cdot v$ by $ v \cdot \frac{u}{v}$).
\end{corollary}
\begin{remark}\label{remark 2}
	Although the Theorem \ref{theorem 1} establishes that $ D(u) = D(\frac{u}{v} \cdot v) = D( v \cdot \frac{u}{v}) $, this should not be interpreted as equality of the underlying ultrafilters. The divisor set captures only the divisibility profile of an ultrafilter and distinct ultrafilters may still possess the same collection of divisors \cite{3}.
\end{remark}
In the arithmetic of natural numbers, divisibility is naturally compatible with quotient formation: whenever $n\mid m$, the quotient $ \frac{m}{n}$ itself remains a divisor of $m$. Thus the quotient preserves an essential part of the divisibility structure of the original number. It is therefore natural to ask whether an analogous phenomenon persists in case of $ \beta \mathbb{N}$. More precisely, if $ v \mid^{s} u$  for ultrafilters $u,v\in \beta\mathbb{N} $, does the quotient $ \frac{u}{v}$ divide  $u$? The answer to this question is given by the following theorem. In particular, as a consequence, one obtains a characterization of self-divisible ultrafilters through the quotient structure.

\begin{theorem}\label{theorem 3}
	Let $ u,v \in \beta \mathbb{N} $ be such that  $ v\mid^{s}u $. Then $u$ is self-divisible if and only if $\frac{u}{v} \mid^{s} u $.
	\begin{proof} If $ u$ is principal then it is trivial. So throughout the theorem  $u$ is assumed to be non-principal.
		\\ $ u $ is self-divisible and $ v \mid^{s} u $ implies $ D(u) \in u $ and $  D(u) \in v $ respectively. Let $ x \in D(u)$. Then $ x\mathbb{N} \cap  D(u)  \in u$. Choose $ y \in x\mathbb{N} \cap  D(u)$ and write $ y = xk $ for some $ k \in \mathbb{N}$. As $ y \in D(u)$, it follows that $ k \in D(u)$ and thus $ y \in xD(u)$. Therefore $ x\mathbb{N} \cap  D(u) \subset  xD(u)$ which further implies that $ x\in \{k\in \mathbb{N}: kD(u) \in u\} $. So,  $ D(u) \subset \{k\in \mathbb{N}: kD(u) \in u\} $ and since $ D(u) \in v $, it implies that the set  $\{k\in \mathbb{N}: kD(u) \in u\} \in v $. Hence $ D(u) \in \frac{u}{v} $ and the result follows.   
		\\[0.3cm]Conversely, $ \frac{u}{v} \mid^{s} u \implies \{x\in \mathbb{N}: xD(u) \in u\} \in v$. Assume that $ u$ is not self-divisible. Pick any $ x(\neq 1) \in \mathbb{N}$ such that $ xD(u) \in u$. Then $D(u)$ is infinite  and consequently $ u$ has infinitely many divisors in $ \mathbb{N}$. Since $u$ is not self-divisible, $ \varphi_u^{-1}(\{0\}) \cup \varphi_u^{-1}(\mathbb{N}) $ must be infinite \cite{3}. As $x \in \mathbb{N}$, it has a finite number of divisors in $P$ so does in $\varphi_u^{-1}(\mathbb{N}) $, if exists. Now there are two cases to consider here.
		\\[0.2 cm] Case (i) : $x$ has divisor(s) in $\varphi_u^{-1}(\mathbb{N}) $.
		\\Let $ p_1, p_2,...,p_k $ be the divisors of $x$ in the set $\varphi_u^{-1}(\mathbb{N})$. Then  $(\bigcap_{i=1}^{k} p_i^{\varphi_u(p_i)} \mathbb{N}) \cap (\bigcap_{i=1}^{k} (p_i^{\varphi_u(p_i)+1} \mathbb{N})^{c})  \cap (\bigcup_{p \in P} p^{\varphi_u(p)+1}\mathbb{N}) \cap \ xD(u) \in u$,  where $  (\bigcup_{p \in P} p^{\varphi_u(p)+1}\mathbb{N}) = (D(u))^c$\cite{3}. Let $ y \in (\bigcap_{i=1}^{k} p_i^{\varphi_u(p_i)} \mathbb{N}) \cap (\bigcap_{i=1}^{k} (p_i^{\varphi_u(p_i)+1} \mathbb{N})^{c})  \cap (\bigcup_{p \in P} p^{\varphi_u(p)+1}\mathbb{N}) \cap xD(u)$. 
		Then $ p_i^{\varphi_u(p_i)} \mid y$ but $p_i^{\varphi_u(p_i)+1} \nmid y$, where  $y = xm $ such that $m \in D(u)$ and $ i \in \{1,2,...k\}$. Moreover $\exists p \in P $ such that $ p^{\varphi_u(p)+1} \mid y $ and $ p \neq p_i, i \in \{1,2,...,k\}$. There are two choices -  either $\varphi_u(p)  = 0$ or $\varphi_u(p) \in \mathbb{N}$.\\ First, suppose  $\varphi_u(p)  = 0$. Then  $ p \mid y$. Since $ m \in D(u)$, it follows that $ m \notin \{\bigcup_{p \in P} p^{\varphi_u(p)+1}\mathbb{N}\}$ and consequently $p \mid x$. Now  $xD(u) \in u$ implies $x\mathbb{N} \in u $ and since $p \mid x $, one has $p\mathbb{N} \in u $ which is not possible because $ \varphi_u(p) = 0$.
		\\ Now come to the second choice, $\varphi_u(p) \in \mathbb{N}$. Then $ p^{\varphi_u(p)+1} \mid y $, whereas $ m \in D(u)$ implies that $  p^{\varphi_u(p)+1} \nmid m$. Therefore it follows that $ p \mid x $. But since $ p \in \varphi_u^{-1}(\mathbb{N}), \ p = p_i $ for some $ i \in \{1,2,...,k\}$, thus arriving at a contradiction. 
		\\Case(ii) : $x$ has no divisor in $\varphi_u^{-1}(\mathbb{N}) $.
		\\Pick any $ y \in (\bigcup_{p \in P} p^{\varphi_u(p)+1}\mathbb{N}) \cap xD(u)$. Then there exists $p \in P $ such that $p^{\varphi_u(p)+1} \mid y $. But since $x$ has no prime divisor in $\varphi_u^{-1}(\mathbb{N})$, $ \varphi_u(p) = 0$. By the previous case, proceeding forward leads to the contradiction.
		\\Therefore, in either case, we arrive at the conclusion that $ \phi \in u$, which is absurd. Hence $u$ must be self-divisible.
	\end{proof}
\end{theorem}	
\begin{lemma}\label{lemma 4}
	Let $u,v,w \in \beta \mathbb{N}$. Then $ w \widetilde{\mid} v, \ v  \mid^{s} u \implies w  \mid^{s} u$.
	\begin{proof}
		Take any $ A \in w$. Then $\bigcup_{n \in A} n \mathbb{N} \in v$. Since $ D(u) \in v$, it follows that  $(\bigcup_{n \in A} n \mathbb{N}) \cap D(u) \in v$. Choose $ x \in (\bigcup_{n \in A} n \mathbb{N}) \cap D(u) $. Then $ x = nk  \in D(u)$ where $n \in A, \ k \in \mathbb{N}$ and consequently $ n \in D(u) \cap A$. Therefore $ D(u)$ intersects every member of $w$. Hence $D(u) \in w $.  
	\end{proof} 
	
\end{lemma}	
\begin{corollary}\label{corollary 6}
	For $u,v \in \beta \mathbb{N}$ with $ v \mid^{s} u$, self- divisibility of $u$ is necessary condition for $u$ to be strongly divisible by $\frac{u}{v} \cdot v$ and $v \cdot \frac{u}{v} $. 
	\begin{proof}
		By Lemma \ref{lemma 4}, $ \frac{u}{v} \cdot v \mid^{s} u $ implies $ \frac{u}{v} \mid^{s} u $ from which it follows that $ u$ is self - divisible. The proof is analogous for $ v \cdot \frac{u}{v} \mid^{s} u $.
	\end{proof}
\end{corollary}

\begin{proposition}\label{proposition 2}
	Let $ u \in \beta \mathbb{N}$ and $m,n \in \mathbb{N}$ such that $ n \mid^{s} u$ and $ m \mid n$. Then $ u $ is self-divisible iff $ \frac{u}{n} \mid^{s} \frac{u}{m}$.
	\begin{proof}  $ n \mid^{s} u $ and $u$ is self-divisible $\implies n\mathbb{N} \cap D(u) \in u $. Pick any $ x \in  n\mathbb{N} \cap D(u)$. Then $ x = kn $, where $ k \in \mathbb{N}$ and $ kn\mathbb{N} \in u$. Now, since $ m \mid n, \ km\mathbb{N} \in u  $ and therefore $ k \in D(\frac{u}{m})$. Thus $ x \in nD(\frac{u}{m}) $ which implies that $  n\mathbb{N} \cap D(u) \subset nD(\frac{u}{m}) $ .  Therefore $ nD(\frac{u}{m}) \in u $, and hence the result follows. 
		\\ Conversely $ \frac{u}{n} \mid^{s} \frac{u}{m} \implies  \frac{u}{n} \mid^{s} u $. From Theorem \ref{theorem 3}, it follows that $u$ is self-divisible.
	\end{proof}
	
\end{proposition}
\begin{remark}\label{remark 3}
	For any $ u \in \beta \mathbb{N} $ and any $ n \in \mathbb{N}$ such that $ n \mid^{s} u, u $ is self-divisible iff $\frac{u}{n}$ is divisible. But the Example \ref{example 1} witnesses that replacement of principal ultrafilter $n$ by non-principal one replaces the conclusion by negation.
	
\end{remark}
\begin{example}\label{example 1}
	
	Let $ S$ denote the set $\{x \in \mathbb{N}: x$ is square-free, that is, $x$ is not divisible by any square of a natural number except $1 \}$. Consider the collection $\mathcal{F} = \{p\mathbb{N} \cap (p^2\mathbb{N})^c: p  \text{ is prime}\} \cup \{S\}$. This collection satisfies FIP and every finite intersection is infinite. So there exists $u \in \beta\mathbb{N} \setminus \mathbb{N}$ such that  $\mathcal{F} \subset u$.
	Then $ D(u)= S $ and since $ S \in u $,  $u$ is self-divisible. But $ P \cap  D(\frac{u}{u}) = \phi$ which implies that $ D(\frac{u}{u}) = \{1\}$. Hence $\frac{u}{u}$ is not self-divisible. 
	
\end{example}
\begin{proposition}\label{proposition 3}
	Let  $ u \in \beta \mathbb{N} $ be such that $u$ is self-divisible as well as multiplicative idempotent. Then for any $v \in  \beta \mathbb{N}, v \mid^{s} \frac{u}{u} \iff v \mid^{s} u $.
	\begin{proof} 
		The forward implication follows directly from Fact \ref{fact 5}.
		\\ Conversely, let $ x \in D(u)$. Take any $y \in D(u)$. This implies $xy\mathbb{N} \in u.u = u$. Thus $ y \in \{k \in \mathbb{N} : kx\mathbb{N} \in u\}$. So, $ D(u) \subset \{k \in \mathbb{N} : kx\mathbb{N} \in u \}$ which means $ x \in D(\frac{u}{u})$. Therefore $ D(u) \subset D(\frac{u}{u})$ and since $ D(u) \in v$, it further implies that $ D(\frac{u}{u}) \in v$. Hence $ v \mid^{s} \frac{u}{u}$. 
		
		\end{proof}	
\end{proposition}
\begin{remark}\label{remark 4}
	As an immediate consequence of Proposition \ref{proposition 3} and Theorem \ref{theorem 3}, it follows that whenever $u$ is a multiplicative idempotent self-divisible ultrafilter, the quotient $\frac{u}{u}$ is itself self-divisible.
\end{remark}
\begin{definition} \label{definition 13} \cite{1}
	Let $u,v\in\beta\mathbb{N}$. Define $E_1(u,v) =\left\{A\subseteq\mathbb{N}:\left\{m\in\mathbb{N}:\left\{n\in\mathbb{N}:m^{n}\in A\right\}\in v\right\}\in u \right\}$,
	and $ E_2(u, v) = \left\{A\subseteq\mathbb{N}:\left\{m\in\mathbb{N}:\left\{n\in\mathbb{N}:n^{m}\in A\right\}\in v\right\}\in u \right\}$.

\end{definition}
In this article, the definition $ E_1(u,v)$ is adopted and is denoted by $ u ^{v}$.
\begin{theorem}\label{theorem 6}
	Let $u \in \beta \mathbb{N} \setminus \mathbb{N} $ be self-divisible. Then for any prime $ p \in \mathbb{N}, \\  p \mid^{s} \frac{u}{u} \iff  p^{k} \mid^{s} u, \forall k \in \mathbb{N}$, i.e, $ \varphi_u(p) = \omega$.
	\begin{proof}
		$ p \mid^{s} \frac{u}{u} \implies \{x \in \mathbb{N} : xp\mathbb{N} \in u\} \in u$. If possible, let $ m \in \mathbb{N}$ be such that $ p^{m} \mid^{s} u$ but $p^{m+1} \nmid^{s} u$. Then, $ p^{m}\mathbb{N} \cap  \{x \in \mathbb{N} : xp\mathbb{N} \in u\} \in u$. Pick $ y \in  p^{m}\mathbb{N} \cap  \{x \in \mathbb{N} : xp\mathbb{N} \in u\}$. As a result $  yp\mathbb{N} \subset p^{m+1}\mathbb{N}$ and it follows that $ p^{m+1}\mathbb{N} \in u$, a contradiction.
		\\[0.2cm]Converse Part : Suppose $ \mathbb{N}\setminus p\mathbb{N} \in \frac{u}{u}$. This implies $ \{x \in \mathbb{N} : x(\mathbb{N}\setminus p\mathbb{N}) \in u\} \in u$. Since $x(\mathbb{N}\setminus p\mathbb{N}) \subset \mathbb{N}\setminus xp\mathbb{N}, \{x \in \mathbb{N} : \mathbb{N}\setminus xp\mathbb{N} \in u\} \in u $. Take any $ m \in \mathbb{N}$. Then $ p^{m}\mathbb{N} \in u$. Let $ y \in (p^{m}\mathbb{N}) \cap \{x \in \mathbb{N} : \mathbb{N}\setminus xp\mathbb{N} \in u\}  \cap D(u) \in u$. Then $y = p^{m}t \in D(u)$, for some $ t \in \mathbb{N}$ and consequently $t \in D(u)$. Now two cases arise.
		\\First, consider the case $ p \nmid t$. In this case $p^{m+1}\mathbb{N} \cap t\mathbb{N} = p^{m+1}t\mathbb{N} = yp\mathbb{N} \in u$, a contradiction.
		Now consider the case $ p \mid t$. There exists $ a,b \in \mathbb{N}$ such that $ t = p^{a}b$ and $ p \nmid b $. It follows that 
		$ b \in D(u)$. Therefore, $ b\mathbb{N} \cap p^{m+a+1}\mathbb{N} = p^{m+a+1}b\mathbb{N} = yp\mathbb{N} \in u $ and this leads to a contradiction. Hence $ p\mathbb{N} \in \frac{u}{u}$.
	\end{proof}
\end{theorem}
\begin{corollary} \label{corollary 7}
	Let $u \in \beta \mathbb{N} \setminus \mathbb{N} $ be self-divisible and $v$ be a prime ultrafilter. Then $ v \mid^{s} \frac{u}{u}$ if and only if $ v^{w} \mid^{s} u \ \forall w \in \beta \mathbb{N}$, 
	\begin{proof}
		Take any ultrafilter $w$. Now, $ v \mid^{s} \frac{u}{u} \implies D(\frac{u}{u}) \in v$. Since $v$ is a prime ultrafilter, $ P \cap D(\frac{u}{u}) \in v$. Let $ p \in P \cap D(\frac{u}{u}) \in v$. By Theorem \ref{theorem 6}, $ \varphi_u(p) = \omega$. Therefore the set $ \{ n \in \mathbb{N} : p^{n} \in D(u)\} = \mathbb{N} $ and it further implies that $P \cap D(\frac{u}{u}) \subset \{m \in \mathbb{N} : \{n \in \mathbb{N} : m^{n} \in D(u)\} \in w\} $. Hence $ D(u) \in v^{w}$.
		\\Conversely, take any $ w \in \beta \mathbb{N} \setminus \mathbb{N}$. Then $ \{m \in \mathbb{N} : \{n \in \mathbb{N} : m^{n} \in D(u)\} \in w\} \cap P \in v$. Now, choose a prime $p$ such that $ \{ n \in \mathbb{N} : p^{n} \in D(u)\} \in w$. Since $w \in \beta \mathbb{N} \setminus \mathbb{N}$, it follows that  the set $\{ n \in \mathbb{N} : p^{n} \in D(u)\}$ is infinite which implies that $ \varphi_u(p) = \omega$, that is, $ p \mid^{s} \frac{u}{u}$. Therefore $\{m \in \mathbb{N} : \{n \in \mathbb{N} : m^{n} \in D(u)\} \in w\} \cap P \subset \{x \in P : x \mid^{s} \frac{u}{u}\} \subset D(\frac{u}{u})$. Hence the result follows.      
	\end{proof}
\end{corollary}
\begin{remark}\label{remark 5}
	In the converse part of the proof of Corollary \ref{corollary 7}, the existence of one non-principal ultrafilter is sufficient to get the required result.
\end{remark}
\begin{corollary}\label{corollary 8}
	For any $ u \in \beta \mathbb{N}$, $ \varphi_u $ is cofinite if and only if  $ \frac{u}{u} $ is strongly divisible by every non-principal prime ultrafilter.
	
	\begin{proof}
		From Theorem \ref{theorem 6}, it follows that $ P \cap D(\frac{u}{u}) = \varphi_u^{-1}(\{\omega\})$. So when $ \varphi_u $ is cofinite, $ P \cap D(\frac{u}{u})$ is cofinite subset of $P$. Hence $D(\frac{u}{u}) \in v$ for every non-principal prime ultrafilter $v$.
		\\ For the converse part, suppose $ \varphi_u$ is not cofinite. Then the set $ P \setminus \varphi_u^{-1}(\{\omega\}) $ is infinite, so there exists a non-principal prime ultrafilter $v$ such that $ P \setminus \varphi_u^{-1}(\{\omega\}) \in v$. But $ \varphi_u^{-1}(\{\omega\}) \in v $ which is a contradiction.
	\end{proof}
\end{corollary}
\begin{corollary}\label{corollary 9}
	If $ u \in \beta \mathbb{N}\setminus \mathbb{N}$ is self-divisible as well as  multiplicative idempotent then : (i) for any prime $ p \in \mathbb{N},p \mid^{s} u \iff  p^{k} \mid^{s} u, \forall k \in \mathbb{N}$, (ii) for any prime ultrafilter $v$, $ v \mid^{s} u \iff v^{w} \mid^{s} u \ \forall w \in \beta \mathbb{N}$.
	\begin{proof}
		Using Proposition \ref{proposition 3} and Theorem \ref{theorem 6}, the result follows.
	\end{proof}
\end{corollary}
\begin{remark}\label{remark 6}
	Theorem \ref{theorem 6} may not hold on replacing $u$ by any other divisor of $u$. Consider the ultrafilter $u$ taken in Example \ref{example 1}. Then $u$ is self-divisible and $u$ is strongly divisible by every prime ultrafilter. Now for any prime ultrafilter $v$, $ P \subset D(\frac{u}{v}) $ but $ \varphi_u(p) = 1 \ \forall p \in P$.	
\end{remark}
In Theorem \ref{theorem 6}, on replacing $p$ by any non-principal prime ultrafilter it is observed that the forward implication holds but the converse may not hold in general as illustrated by Example 2.
\begin{example}\label{example 2}
	Consider the collection $ \mathcal{F} = \{p^{p}\mathbb{N} \cap (p^{p+1}\mathbb{N})^{c} : p \in P\}$ and any prime ultrafilter $v$. $ \mathcal{F}$ satisfies FIP, so there exists $u \in \beta \mathbb{N}$ such that $ \mathcal{F} \subset u$. Fix any $n \in \mathbb{N}$. The set $\{p \in P : p > n\} \subset \{m \in \mathbb{N} : m^{n} \in D(u)\}$. As $\{p \in P : p > n\}$ is cofinite in $P, \{m \in \mathbb{N} : m^{n} \in D(u)\} \in v$. Since this holds for any $n, \ v^{n} \mid^{s} u \ \forall n \in \mathbb{N}$. But $ \varphi_u(p) = p \ \forall p \in P$ implies that $ P \cap D(\frac{u}{u}) = \varphi_u^{-1}({\omega}) = \phi$. Therefore $ v \nmid^{s} \frac{u}{u}$.
\end{example}
\begin{proposition}\label{proposition 4}
	Let $u \in \beta \mathbb{N}$ be a multiplicative idempotent. Then for any $v,w \in \beta \mathbb{N} $ such that $ v \mid^{s} u $ and $ w \mid^{s} v, u $ is self-divisible iff $ \frac{u}{v} \mid^{s} \frac{u}{w}.$
	\begin{proof} $u$ is self-divisible and $ v \mid^{s} u \implies D(u) \in u$ and $D(u) \in v $ . By transitivity of $\mid^{s}, w \mid^{s} u$ which means  $D(u) \in w$.  Take any $ x \in D(u) $.  Then $ x\mathbb{N} \cap D(u) \in u $. Now pick any $ y \in x\mathbb{N} \cap D(u)$. Write $ y = xk $, where $ k \in \mathbb{N} $ and $ y\mathbb{N} \in u$. It follows that $ k\mathbb{N} \in u $. For any $ n \in D(u) , nk\mathbb{N} \in u (\because u\cdot u =u)$, which implies $ n \in \{m \in \mathbb{N}: m(k\mathbb{N}) \in u\} $. Thus $ D(u) \subset  \{m \in \mathbb{N}: m(k\mathbb{N}) \in u\} $ and since $ D(u) \in w $, it follows that $ \{m \in \mathbb{N}: m(k\mathbb{N}) \in u\} \in w $, that is, $k \in D(\frac{u}{w})$. Consequently $  x\mathbb{N} \cap D(u)  \subset x D(\frac{u}{w})$ and thus $ xD(\frac{u}{w}) \in u $. It further implies that $ x \in \{m \in \mathbb{N}: mD(\frac{u}{w}) \in u\} $ . Therefore $ D(u) \subset  \{m \in \mathbb{N}: mD(\frac{u}{w}) \in u\} $. As \ $ D(u) \in v $  it implies  $ D(\frac{u}{w}) \in \frac{u}{v}$. Hence the result follows.
		\\Conversely $ \frac{u}{v} \mid^{s} \frac{u}{w} \implies D(\frac{u}{w}) \in \frac{u}{v} $.  Since $ D(\frac{u}{w}) \subset D(u) $, it follows that $ D(u) \in \frac{u}{v} $, that is, $ \frac{u}{v} \mid^{s} u$. Hence from Theorem \ref{theorem 3}, self-divisibility of $u$ follows.
	\end{proof}
	
\end{proposition}

\begin{remark}\label{remark 7}For self-divisible ultrafilters $u,v $ such that $ v \mid^{s} u, \ \frac{u}{v} $ may not be self-divisible, see Example 1. By considering $ u $ as a multiplicative idempotent and self-divisible, it is immediate from Proposition 4 that for any self-divisible ultrafilter $v $ such that  $ v \mid^{s} u,\ \frac{u}{v} $  is self-divisible.
\end{remark}
The preceding results reveal several structural relations among the classes of self-divisible ultrafilters, additive idempotents, multiplicative idempotents and the class $MAX$. For convenience, we summarize these relations schematically in the following diagram (Figure 1), where $ SD $ \cite{3} denotes the class of self-divisible non-principal ultrafilters, $ \mathcal{A}$ denotes the class of additive idempotents, $\mathcal{M}$ denotes the class of multiplicative idempotents and $ MAX $ is the class of ultrafilters divisible by all the ultrafilters (see Definition 6).
\vspace{-0.7 cm}
\begin{center}
	\begin{tikzpicture}[scale=1.05]
		
		
		\draw[blue, thick] (0,0) ellipse (6 and 4);
		
		\node[blue] at (0,3.1)
		{\fontsize{28}{30}\selectfont $\mathbf{SD}$};

		
		\draw[green!60!black, thick]
		(-1.8,-0.7) ellipse (3.7 and 2.4);
		
		\node[green!50!black] at (-2,1.1)
		{\fontsize{28}{30}\selectfont $\mathbf{MAX}$};

		
		\draw[red, thick]
		(-3.2,-1.1) ellipse (1.6 and 1.25);
		
		\node[red] at (-3.2,-0.9)
		{\fontsize{28}{30}\selectfont $\mathcal{A}$};

		
		\draw[orange!90!black, thick]
		(3.6,-1.3) ellipse (3.4 and 2.5);
		
		\node[orange!80!black] at (4.1,-1.0)
		{\fontsize{28}{30}\selectfont $\mathcal{M}$};

		
		\node[circle,draw,fill=white,thick]
		at (-4,-1.0)
		{\Large $I$};
		
		\node[circle,draw,fill=white,thick]
		at (-0.3,-0.1)
		{\Large $II$};
		
		\node[circle,draw,fill=white,thick]
		at (1.0,-1.15)
		{\large $III$};

		\node[circle,draw,fill=white,thick]
		at (3.2,2.0)
		{\large $IV$};
		
		\node[circle,draw,fill=white,thick]
		at (3,-1.15)
		{\Large $V$};
		
		\node[circle,draw,fill=white,thick]
		at (5.6,-2.6)
		{\large $VI$};
	\end{tikzpicture}
\end{center}
\textbf{Figure 1} : Structural relations among the classes of self-divisible ultrafilters, additive idempotents, multiplicative idempotents and the class $MAX$.
\begin{itemize}
	\item Example for I : Any additive idempotent ultrafilter. 
	\item Example for II : Let $ A = \{n! : n \in \mathbb{N}\}$.The collection  $ \mathcal{C} = \{n\mathbb{N}: n \in \mathbb{N}\} \cup \{A\}$ has the Finite Intersection Property, so there exists $ u \in \beta \mathbb{N}\setminus \mathbb{N}$ such that $ \mathcal{C} \subset u$.
	Then $ u \in MAX $. $ A $ is not an $IP$-set (neither additive nor multiplicative) and therefore it is immediate from Fact \ref{fact 3}  that $u$ is not an idempotent.
	\item Example for III : $ \mathcal{A} \subset MAX $ and $ MAX $ is closed in $\beta \mathbb{N}$ implies $ \overline{\mathcal{A}} \subset MAX$. Since  $ \overline{\mathcal{A}}$ forms a compact Hausdorff subsemigroup of $ (\beta \mathbb{N}, \cdot)$ \cite{4}, there exists a multiplicative idempotent $ u \in \overline{\mathcal{A}} \subset MAX $. 
	\item Example for IV : See Example \ref{example 1}.
	\item Example for V : 	Consider the set $M = \{2^{n}: n \in \mathbb{N} \}$. $\widehat{M}$ forms a subsemigroup of $\beta \mathbb{N}$, also $\widehat{M}$ is compact and Hausdorff. Therefore there exists a multiplicative idempotent $u \in \widehat{M}$. For each $x \in M$, write $ x = 2^{k}$ where $ k \in \mathbb{N}$. Then $M\setminus\{2,...,2^{k-1}\} \subset 2^{k}\mathbb{N} = x\mathbb{N}$ and this implies $ x \in D(u)$. Therefore $ M \subset D(u)$ and thus $D(u) \in  u$, that is, $u$ is self-divisible. Since $3\mathbb{N} \cap M = \phi, \ 3\mathbb{N} \notin u $. Hence $ u \notin MAX$.
	\item Example for VI : See Proposition 5.9 in \cite{3}.
\end{itemize}

\begin{definition}\label{definition 14}
	\cite[1.5. Definition]{2}	$A \subset \mathbb{N} $ is said to be a multiplicative $\Delta$-set if there exists a strictly increasing sequence $(x_n)$ in $\mathbb{N}$
	with the property $ x_n \mid x_m $ whenever $ m > n$ , such that $\{\frac{x_m}{x_n}:m>n\} \subset A$.
\end{definition}

Note: In the definition of a $\Delta$-set, we require the underlying sequence to be strictly increasing, although no such assumption is imposed in \cite{2}. This additional condition is introduced to exclude trivial examples arising from constant sequences. Indeed, without this condition, every subset of $\mathbb{N}$ containing $1$ would automatically become a $\Delta$-set by considering the constant sequence. Therefore the extra condition on sequence ensures that the notion reflects genuine combinatorial largeness.
\begin{definition}\label{definition 15}
	\cite[1.8. Definition]{2}	$A \subset \mathbb{N} $ is said to be a multiplicative $\Delta^*$-set if $A\cap B \neq \phi$  for every $\Delta$-set $B$ in $ \mathbb{N} $.	
\end{definition}
The notion of self-divisibility is closely related to multiplicative combinatorial largeness. The following theorem characterizes  multiplicative $ \Delta$ -sets in terms of quotient of self-divisible ultrafilter.

\begin{theorem}\label{theorem 10}
	Let $A \subset \mathbb{N}$. Then $A $ is a $\Delta$-set if and only if $ A \in \frac{u}{u} $ for some self-divisible non-principal ultrafilter $u$.
	\begin{proof}By assuming $A$ is a $\Delta$-set, there exists a strictly increasing sequence $(x_n)$ in $\mathbb{N}$ with the property $ x_n \mid x_m $ whenever $ m > n$ such that $\{\frac{x_m}{x_n}:m>n\} \subset A$.
		\\ For each $n \in \mathbb{N}$, define $ B_n = \{x_m: m>n\}$. The family  $\{B_n: n \in \mathbb{N} \}$ has the finite intersection property, so there exists $ u \in \beta \mathbb{N}\setminus \mathbb{N}$ such that $\{B_n: n \in \mathbb{N} \} \subset u$. Consequently the set $ B = \{x_n : n \in \mathbb{N} \} \in u$.
		Take any $ x_n \in B $ for some $n \in \mathbb{N}$. Then $ B_n \subset$ $x_n$$A$. Since $B_n \in u$, $x_n$$A$  $\in u$. Thus $B \subset \{x\in \mathbb{N}: xA \in u \}$ and this implies that $\{x\in \mathbb{N}: xA \in u \} \in u$. Hence $A \in \frac{u}{u}$.
		\\ Conversely, let $ X = \{x\in \mathbb{N}: xA \in u \}$. Then $ X \in u$.  Choose  $x_1 \in X$. As $x_1A$ $\in u$, $X \cap $ $x_1A$ $\in u $. Choose $x_2 > x_1 $ such that $x_2 \in X \cap $ $x_1A$. Proceeding inductively, we obtain a strictly increasing sequence $(x_n)$ in $\mathbb{N}$ such that $ x_n \mid x_m $ whenever $ m > n$ and $x_{n+1} \in X \cap x_nA \cap...\cap x_1A  $ for every $ n \in \mathbb{N}$. Therefore $\{\frac{x_m}{x_n}:m>n\} \subset A$ showing that $A$ is a $\Delta$- set .
	\end{proof}
	
\end{theorem}
The technique used in the proof of the Theorem \ref{theorem 10} is quite similar to the one developed by Bergelson and Hindman in \cite{1}.
\begin{corollary}\label{corollary 11}
	Let $A \subset \mathbb{N}$. Then $A $ is a $\Delta^*$-set if and only if $A \in \frac{u}{u} $ for every self-divisible non-principal ultrafilter $u$.
	\begin{proof}If possible, let $ u \in \beta \mathbb{N}\setminus \mathbb{N}$ be such that $u$ is self-divisible and $A \notin \frac{u}{u}$. Then $\mathbb{N}\setminus A \in \frac{u}{u}$. By Theorem \ref{theorem 10}, it follows that $\mathbb{N}\setminus A $ is a $\Delta$-set. But $ A \cap(\mathbb{N}\setminus A)   = \phi$, contradicting the assumption that $ A $ is a $\Delta^{*}$- set\\For the converse, let $B \subset \mathbb{N}$ be any $\Delta$-set. Then there exists a self-divisible non-principal ultrafilter $u$ such that $ B \in \frac{u}{u} $. Since $ A \in \frac{u}{u} $, $ A \cap B $ is non-empty. Hence, $ A$ is a $ \Delta^*$-set.
	\end{proof}	
\end{corollary}

\begin{corollary}\label{corollary 12}
	There is no prime ultrafilter of the form $ \frac{u}{u}$, where $u$ is a self-divisible non-principal ultrafilter.
	\begin{proof} Suppose  there exists a self-divisible ultrafilter $u \in \beta \mathbb{N}\setminus \mathbb{N}$ such that $\frac{u}{u}$ is a prime ultrafilter. Then $P$ is a $\Delta$-set. So there exists a sequence $(x_n)$ in $\mathbb{N}$ with the property $ x_n \mid x_m $ whenever $ m > n$, such that $\{\frac{x_m}{x_n}:m>n\} \subset P$. Let $\frac{x_2}{x_1}=p_1 \in P, \ \frac{x_3}{x_2}=p_2 \in P $. Then $\frac{x_3}{x_1}= \frac{x_2}{x_1} \cdot \frac{x_3}{x_2}=p_1p_2 \notin P$, a contradiction.
	\end{proof}
\end{corollary}

\begin{proposition}\label{proposition 5}
	Let $ v \in \beta \mathbb{N} $ be a prime ultrafilter and let $u_1 , u_2 \in  \beta \mathbb{N} $. Then $ v\mid^{s} u_1\cdot u_2 \iff $ either $ v \mid^{s} u_1 $ or $ v \mid^{s} u_2 $.
	\begin{proof} $ v \mid^{s}  u_1\cdot u_2 \implies D(u_1\cdot u_2) \in v \implies \{x \in \mathbb{N}: x\mathbb{N} \in  u_1\cdot u_2 \} \in v$. Let $ p \in P \in v$ such that $ p \in \{x \in \mathbb{N}: x\mathbb{N} \in  u_1\cdot u_2 \}$. Then $ p\mathbb{N} \in  u_1\cdot u_2$. If neither $ p\mathbb{N} \in u_1 $ nor $ p\mathbb{N} \in u_2$, then $ p\mathbb{N} \in u_1\cdot u_2$ implies $\{x \in \mathbb{N}: p\mathbb{N}/x \in u_2\} \in u_1 $ and $ p\mathbb{N} \notin u_1$ implies $\mathbb{N}\setminus p\mathbb{N} \in u_1 $.
		Pick $ x \in \mathbb{N}\setminus p\mathbb{N} $ such that $ p\mathbb{N}/x  \in u_2$. But since $ x \in \mathbb{N}\setminus p\mathbb{N},\  p\mathbb{N}/x = p\mathbb{N}$. Thus $p\mathbb{N} \in u_2$, contradicting our assumption. Therefore, \ $ p \in D(u_1) \cup D(u_2) $, which shows that $ P \cap D(u_1\cdot u_2) \subset D(u_1) \cup D(u_2) $. Now, $  D(u_1) \cup D(u_2) \in v $ implies either $D(u_1) \in v $ or $D(u_2) \in v $. 
		\\ Reverse implication is straightforward from Fact \ref{fact 2}.
		
	\end{proof}
\end{proposition}
\begin{corollary}\label{corollary 13}
	Let $ v \in \beta \mathbb{N}$ be such that $S = \{ x \in \mathbb{N}: x \text{ is square-free}\} \in v$ and let $ u \in \beta \mathbb{N}$, then $ v \mid^{s} u\cdot u \iff v \mid^{s} u$.
	\begin{proof}$ v \mid^{s} u \cdot u \implies \{ x \in \mathbb{N} : x\mathbb{N} \in u\cdot u\} \in v$. Hence $ S \cap \{ x \in \mathbb{N} : x\mathbb{N} \in u\cdot u\} \in v$. Choose $x \in S \cap \{ x \in \mathbb{N} : x\mathbb{N} \in u\cdot u\} $.Then $x $ is either prime or product of distinct primes. If $ x $ is prime then by Proposition \ref{proposition 5}, $ x \in D(u) $. 
		If $x $ is not a prime then $ x = p_1...p_k $ where $p_1, ..., p_k $ are distinct primes.
		\\ Now, $ x\mathbb{N} = p_1...p_k\mathbb{N} = p_1\mathbb{N} \cap ...\cap p_k\mathbb{N} \in u\cdot u \implies p_i\mathbb{N} \in u\cdot u$, for each $ i \in \{1,...,k\}$ and applying Proposition \ref{proposition 5}, $p_i\mathbb{N} \in u$, for each $ i \in \{1,...,k\}$ which implies $x \in D(u)$. Therefore $ S \cap \{ x \in \mathbb{N} : x\mathbb{N} \in u\cdot u\} \subset D(u)$ and consequently $ D(u) \in v$. Hence $ v \mid^{s} u$.
		\\ Converse part is immediate from Fact \ref{fact 2}.
		
	\end{proof}
\end{corollary}
The class of self-divisible ultrafilters  plays a central role in the quotient framework developed so far in this article. By Theorem \ref{theorem 10}, self-divisibility is intimately linked to the combinatorial structure of multiplicative $\Delta$-sets. Moreover, for every $u\in SD$, the quotient $\frac{u}{u}\in \beta \mathbb{N}$, giving rise to a distinguished family of ultrafilters. Define $\mathcal{D}= \{\frac{u}{u}: u \in SD\}$. As shown in \cite[Proposition 6.3]{3}, $\overline{SD} = \{ u \in \beta \mathbb{N}: \forall A \in u,  \exists X \subset A ( X \text{is an infinite} \mid - chain)\}$, and by \cite{4}, for any property $P$, the set $ \{u \in \beta \mathbb{N} : $ every member of $u$ satisfies $P\}$ is closed in $ \beta \mathbb{N}$. These facts are applied in the next results to characterize the closure of $\mathcal{D}$ and to examine its relationship with the closure of $SD$.
\begin{proposition}\label{proposition 6}
	(i) $\overline{\mathcal{D}} = \{v \in \beta \mathbb{N}:\text{every member of $v$ is a $\Delta$-set}\}$, (ii)  $\overline{\mathcal{D}} \subset \overline{SD}$.
	\begin{proof} (i) By Theorem \ref{theorem 10}, it follows that $\overline{\mathcal{D}} \subset \{v \in \beta \mathbb{N}:\text{every member of $v$ is a $\Delta$-set}\}$. For the reverse inclusion, let $ v \in \beta \mathbb{N} $ be such that every member of $v$ is a $\Delta$-set. Then, for any $A\in v$, $A \in \frac{u}{u} $ for some self-divisible non-principal ultrafilter $u$. Hence $ v \in \overline{\mathcal{D}}$.
		\\[0.2 cm](ii) Let $ u \in \overline{\mathcal{D}}$, and $ A \in u$. Then there exists a strictly increasing sequence $(x_n)$ in $ \mathbb{N}$ with the property $x_n \mid x_m$ whenever $ m > n $, such that $ \{\frac{x_m}{x_n}: m>n\} \subset A$. Define a sequence $ (y_n) $ such that $ y_n = \frac{x_{n+1}}{x_1} $. For any $m,n \in \mathbb{N}$ with $m>n$, we have $ y_m = \frac{x_{m+1}}{x_1} = \frac{x_{m+1}}{x_m} \cdot \frac{x_m}{x_{m-1}}\cdot \cdot \cdot \frac{x_{n+1}}{x_1}$ which implies $ y_n \mid y_m$. Hence, the set $\{y_n: n \in \mathbb{N}\}$ forms an infinite divisibility chain contained in $A$. Therefore, $ u \in \overline{SD}$.
	\end{proof}
	
\end{proposition}
The next Example 3 shows that $\overline{\mathcal{D}}$ is not necessarily equal to $\overline{SD}$.
\begin{example}\label{example 3}
	Consider the set $\mathbb{A}=\{2^{2^{n}}: n \in \mathbb{N}\}$. There exists a non-principal ultrafilter $u$ such that  $\mathbb{A} \in u$.  So, $\{2^n : n \in \mathbb{N}\} \in u $. From Example V, it follows that $u \in  SD$. If possible let, there exists a strictly increasing sequence $(x_n)$ in $\mathbb{N}$ such that $ \{\frac{x_m}{x_n}: m>n\} \subset \mathbb{A} $. Then $ \frac{x_2}{x_1} = 2^{2^{a_1}}$, $\frac{x_3}{x_2} = 2^{2^{a_2}}$, and $\frac{x_4}{x_3} = 2^{2^{a_3}} $ for some $a_1, a_2, a_3 \in \mathbb{N}$. As $ \frac{x_2}{x_1} \cdot \frac{x_3}{x_2} =\frac{x_3}{x_1}\in \mathbb{A}$, it follows that  $a_1 = a_2 $. Similarly,  $a_2 = a_3  $. Now, $ \frac{x_2}{x_1} \cdot  \frac{x_3}{x_2} \cdot  \frac{x_4}{x_3}=\frac{x_4}{x_1}  \in \mathbb{A}$. Therefore, $2^{3\cdot{2^{a_1}} } \in \mathbb{A}$, which is impossible. Hence, $ u \notin \overline{\mathcal{D}} $.
\end{example}
Result analogous to the following theorem was previously obtained by B. \v{S}obot \cite[Corollary 6.3(b)]{7} in the context of $\widetilde{\mid}$ divisibility.	
\begin{theorem}\label{theorem 14}
	A self-map on $\mathbb{N}$ preserving standard divisibility, when extended continuously to $\beta \mathbb{N}$, preserves strong divisibility.
	\begin{proof} Suppose $f:\mathbb{N}\to \mathbb{N}$ is a function satisfying $f(a)\mid f(b)$ whenever $a\mid b$. Then, $ f(x\mathbb{N}) \subset f(x)\mathbb{N} $ for all $x\in \mathbb{N}$. For $u,v\in \beta \mathbb{N}$ with $v\mid^{s} u$, $ D(u) \in v $ and so $f(D(u)) \in \widetilde{f}(v)$. Further, if $x\in D(u)$, then $x\mathbb{N}\in u$, and hence $f(x\mathbb{N}) \in \widetilde{f}(u)$. Therefore, $f(D(u)) \subset D(\widetilde{f}(u))$. This implies $ D(\widetilde{f}(u)) \in \widetilde{f}(v)$. Thus, $\widetilde{f}(v) \mid^{s} \widetilde{f}(u)$.
		
	\end{proof}	
\end{theorem}
\begin{corollary}\label{corollary 15}
	For any multiplicative homomorphism $ f 	: \mathbb{N} \rightarrow \mathbb{N} $ and any  $ u, v \in \beta \mathbb{N}$, $\widetilde{f}(v) \mid^{s} \widetilde{f}(u)$ whenever $ v \mid^{s} u$.
\end{corollary}
The compatibility between quotient structure arising from ultrafilter divisibility and Stone-\v{C}ech extensions of self mappings on $ \mathbb{N}$ leads to a characterization of multiplicative homomorphisms. In particular, the quotient structure is preserved under Stone-\v{C}ech extensions induced by multiplicative homomorphisms, as established in the next theorem.
\begin{theorem}\label{theorem 16}A function $f : \mathbb{N} \to \mathbb{N}$ is a multiplicative homomorphism precisely when $  \widetilde{f}(\frac{u}{v})= \frac{\widetilde{f}(u)}{\widetilde{f}(v)}$ for any $u,v \in \beta \mathbb{N}$ satisfying $ v \mid^{s} u $.
	\begin{proof}
		Let $ f: \mathbb{N} \rightarrow \mathbb{N} $ be a multiplicative homomorphism. Suppose $u,v\in \beta \mathbb{N}$ such that $ v \mid^{s} u $. By Corollary \ref{corollary 15}, $\frac{\widetilde{f}(u)}{\widetilde{f}(v)} \in \beta \mathbb{N}$. Pick  any $ A \in \widetilde{f}(\frac{u}{v}) $. As $f^{-1} (A) \in \frac{u}{v}$, $B=\{x \in \mathbb{N}: xf^{-1} (A) \in u\}\in v$, which follows that $f(B) \in \widetilde{f}(v)$. Since $ f$ is a multiplicative homomorphism, $f(xf^{-1} (A)) \subset f(x)A$ for every $x\in \mathbb{N}$. Moreover, $f(xf^{-1} (A)) \in \widetilde{f}(u)$ for every $x\in B$, and so $ f(x)A \in \widetilde{f}(u)$. Hence, $ f(B) \subset \{k \in \mathbb{N}: kA \in \widetilde{f}(u)\}$. Therefore, $\{k \in \mathbb{N}: kA \in \widetilde{f}(u)\} \in \widetilde{f}(v)$, whence $ A \in \frac{\widetilde{f}(u)}{\widetilde{f}(v)} $ and consequently $\widetilde{f}(\frac{u}{v}) = \frac{\widetilde{f}(u)}{\widetilde{f}(v)}$.
		
		Conversely, $f(\frac{ab}{a}) = \frac{f(ab)}{f(a)}$ for all $ a, b \in \mathbb{N}$, as $ a \mid ab$ and  $\widetilde{f}\vert_\mathbb{N} = f$. Therefore, $ f(ab) = f(a)f(b) $.
	\end{proof}
\end{theorem}
The concept of quotient developed in this work is specific to strong divisibility.  It is natural to investigate whether similar quotient constructions can be formulated for other extensions of the standard divisibility relation on $\mathbb{N}$ to $ \beta \mathbb{N}$, including $\widetilde{\mid}$ - divisibility, left divisibility and right divisibility. Since these divisibility relations possess different algebraic and topological characteristics, the corresponding quotient structures are expected to exhibit significantly different properties which may give sight to a new understanding of divisibility relations on $\beta \mathbb{N}$ and provide new directions for further research.

\end{document}